\documentclass[]{article}
\usepackage{amsmath}
\usepackage{amssymb}
\usepackage{bbm}
\usepackage[margin=2.25cm]{geometry}
\usepackage{hyperref}
\usepackage{tikz}
\usepackage{setspace}

\usetikzlibrary { decorations.pathmorphing, decorations.pathreplacing, decorations.shapes, }

\newtheorem{thm}{Theorem}[section]
\newtheorem{cor}[thm]{Corollary}
\newtheorem{prop}[thm]{Proposition}

\newtheorem{lemma}[thm]{Lemma}

\newtheorem{theorem}[thm]{Theorem}
\newtheorem{remark}[thm]{Remark}
\newtheorem{proposition}[thm]{Proposition}
\newtheorem{definition}[thm]{Definition}

\newenvironment{proof}{{\bf Proof:}}{\hfill$\square$\vskip.5cm}

\newcommand{\R}{\mathbb{R}}
\newcommand{\N}{\mathbb{N}}

\newcommand{\Z}{\mathbb{Z}

}

\renewcommand{\Z}{{\mathbb{Z}}}

\begin{document} 

\title{Two-type thinning with alternative leader rate for random partition structures}
\author{Gracelyn Varghese, Katherine Graham, and Shannon Starr\\
\small
UAB, Math.~Dept., 1402 Tenth Avenue South, Birmingham, AL 35294--1241\\
}
\date{23 November 2026}

\maketitle

\begin{abstract}
A dust-free random partition structure with infinitely many parts  is a random sequence
$(\mathsf{X}_1,\mathsf{X}_2,\dots)$ with $\mathsf{X}_1>\mathsf{X}_2>\dots>0$ such that $\mathsf{X}_1+\mathsf{X}_2+\dots=1$,
almost surely. The Poisson-Dirichlet $\mathrm{PD}(m,0)$ distributions for $0<m<1$ are defined such that
$\mathsf{X}_n = \xi_n^{-1/m}/\sum_{\nu=1}^{\infty} \xi_{\nu}^{-1/m}$, where $0<\xi_1<\xi_2<\dots$ is a  homogeneous Poisson point
process.
It is a model for thinning invariant random partition structures, where every one of the infinitely many parts $\mathsf{X}_1,\mathsf{X}_2,\dots$
are thinned IID at rate 1.
We consider a 2-type thinning where the instantaneous leader is deleted with an alternative rate $1+r\geq 0$.
Then we seek time stationary random partition structures.
When $r =R\in \{0,1,\dots\}$  a solution is obtained by taking $\mathrm{PD}(m,0)$ and performing ``balayage,'' 
where balayage means simply
excavating the first $R$ points.
We find the solutions for all $r \in (-1,\infty)$, which we call
$\mathrm{PD}^{(r)}(m,0)$.
We motivate with an application, and we present graphical representations.

\end{abstract}

\section{Balayage for partition structures}

Let us denote $\mathcal{D}_{\infty}^{+}$ as the set of sequences $(x_1,x_2,\dots)$ such that
$$
x_1\, >\, x_2\, >\, \dots\, >\, 0\, .
$$
Let $\widehat{\mathcal{D}}_{\infty}^{+}$ be the subset of $\mathcal{D}_{\infty}^{+}$ consisting of those sequences
$(x_1,x_2,\dots)$ which satisfy
$$
\sum_{n=1}^{\infty} x_n\, <\, \infty\, .
$$
Of course, the sum is strictly positive.
Let $\overline{\mathcal{D}}_{\infty}^{+}$ be the subset of $\widehat{\mathcal{D}}_{\infty}^{+}$ such that the sum satisfies
$$
\sum_{n=1}^{\infty} x_n\, =\, 1\, .
$$
Then a dust-free random partition structure with infinitely many points is a probability measure on $\overline{\mathcal{D}}_{\infty}^{+}$.
The random element is $\mathsf{X}_1>\mathsf{X}_2>\dots>0$.

\begin{definition}
\label{def:PDm0}
    Let $0<\xi_1<\xi_2<\dots$ be a standard, homogeneous Poisson point process (so that $\xi_1,\xi_2-\xi_1,\xi_3-\xi_2,\dots$ are IID exponential-1 random variables).
    For $q>1$, let $\Phi_q : (0,\infty) \to (0,\infty)$ be defined such that $\Phi_q(x) = x^{-q}$.
    Then define $\mathsf{X}_1,\mathsf{X}_2,\dots$ by the formula
    $$
        \mathsf{X}_n\, =\, \frac{\Phi_q(\xi_n)}{\sum_{\nu = 1}^{\infty} \Phi_q(\xi_{\nu})}\, ,
    $$
    for $n \in \N = \{1,2,\dots\}$. The law of the random partition structures $(\mathsf{X}_1,\mathsf{X}_2,\dots)$ is called Poisson-Dirichlet
    $\mathrm{PD}(m,0)$ for $m=1/q$.
\end{definition}

\begin{remark}
    It is easy to see that $\Phi_q(\xi_1)>\Phi_q(\xi_2)>\dots$ is a Poisson point process on $(0,\infty)$ with intensity measure $mx^{-m-1}\, dx$ for $m=1/q$.
    Moreover, $\sum_{n=1}^{\infty} \Phi_q(\xi_n)$ is a positive one-sided $m$-stable random variable with 
    $\mathbb{E}[\sum_{n=1}^{\infty} \Phi_q(\xi_n)]=\exp(-\lambda^m \Gamma(1-m))$.
    In particular, it is finite, a.s. So the description above of the random partition structure $\mathrm{PD}(m,0)$ is well-defined.
\end{remark}

An excellent reference is Pitman and Yor \cite{PitmanYor}.

Given a  partition structure with infinitely many points $(x_1,x_2,\dots) \in \overline{\mathcal{D}}^{+}_{\infty}$,
let us say the ``$R$-point balayage,'' for $R \in \{0,1,\dots\}$, is $(y_1,y_2,\dots) \in \overline{\mathcal{D}}^{+}_{\infty}$
where
$$
y_n\, =\, \frac{x_{n+R}}{\sum_{\nu=1}^{\infty} x_{\nu+R}}\, .
$$
If $\textsc{Law}$ is the law of a random partition structure $(\mathsf{X}_1,\mathsf{X}_2,\dots)$
let us write $\mathcal{B}^{R}(\textsc{Law})$ for the law of $(\mathsf{Y}_1,\mathsf{Y}_2,\dots)$
where
$$
\mathsf{Y}_n\, =\, \frac{\mathsf{X}_{n+R}}{\sum_{\nu=1}^{\infty} \mathsf{X}_{\nu+R}}\, .
$$
Now we will define a new random partition structure related to the Poisson-Dirichlet 
$\mathrm{PD}(m,0)$ process.

\begin{definition}
    For $r \in (-1,\infty)$, let  $\xi_1^{(r)}$ be an independent Gamma-$(1+r,1)$ random variable.
    Then, conditional on the value of $\xi^{(r)}_1$, let $\xi^{(r)}_2<\xi^{(r)}_3<\dots$ be distributed
    as a homogeneous Poisson point process of intensity 1 on $(\xi^{(r)}_1,\infty)$.
    Finally, take
    $$
    \mathsf{X}_n^{(r)}\, =\, \frac{\Phi_q(\xi_n^{(r)})}
    {\sum_{\nu=1}^{\infty}\Phi_q(\xi_\nu^{(r)})}\, .
    $$
    We denote the distribution as $\mathrm{PD}^{(r)}(m,0)$.
\end{definition}

If $r \in (0,\infty)$ another alternative to this definition exists which we will use more
frequently.
Let $\mathsf{Z}^{(r)}$ be an independent Gamma-$(r,1)$ random variable
and let $\xi_1<\xi_2<\dots$ be a homogeneous Poisson point process of rate $1$ in $(0,\infty)$.
Then let $\xi^{(r)}_n = \mathsf{Z}^{(r)}+\xi_n$.
Note that for $r=0$, we may simply take $\mathsf{Z}^{(0)}=0$, a.s.

\begin{lemma}
For each $R \in \{0,1,\dots\}$, we have $\mathcal{B}^{R}(\mathrm{PD}(m,0))$
equals $\mathrm{PD}^{(R)}(m,0)$.
\end{lemma}
\begin{proof}
    The distribution of $\xi_R$ is equal in distribution to $\mathsf{Z}^{(R)}$ because the sum of $R$ independent exponential-1 random variables is Gamma-$(R,1)$. And the gaps of both processes are
    independent of the lead point, and are distributed as IID Exponential-1 random variables.
\end{proof}

\section{2-type thinning with alternative leader rate 1+r}

\subsection{Motivation: market leadership and regulatory scrutiny}
\label{subsec:motivation}
Imagine an emerging sector with a large number of competing platforms or firms $\mu_1, \mu_2, \dots$, ordered by decreasing market share $X_1 > X_2 > \cdots > 0$. For simplicity, we model the sector with infinitely many firms satisfying $\sum_{n=1}^\infty X_n = 1$. Each firm initially has status $A$, denoting \textsc{active competitor}, and eventually transitions to status $C$, denoting \textsc{exited}. For $r > 0$, a firm may also pass through an intermediate status $B$, denoting \textsc{restricted} (for instance, operating under regulatory limits or confined to a legacy niche), before exiting to $C$.

In a symmetric baseline market, routine competitive attrition, acquisition, or obsolescence causes every firm to exit to state $C$ independently at unit rate, whether it is instantaneously in state $A$ or state $B$. We call this background attrition \emph{type-1 thinning}.

In practice, however, the instantaneous market leader rarely experiences the same exit hazard as smaller competitors. At time $t \geq 0$, define $\pi^{(r)}(\cdot, t) : \mathbb{N} \to \mathbb{N}$ by $\pi^{(r)}(1,t) < \pi^{(r)}(2,t) < \cdots$ and
\[
  \{\pi^{(r)}(1,t), \pi^{(r)}(2,t), \dots\} = \{n \in \mathbb{N} : \mu_n \text{ is in state } A \text{ at time } t\},
\]
so that $\mu_{\pi^{(r)}(1,t)}$ is the leading active firm at time $t$. While smaller firms often operate below the threshold of public or governmental attention, reaching the top rank alters the leader's hazard rate to $1 + r > 0$, giving two distinct regimes for $r > -1$:
\begin{itemize}
  \item \textbf{Targeted regulatory or competitive pressure ($r > 0$):} The visibility and scale of the dominant platform attract heightened antitrust scrutiny, targeted regulatory intervention, or concentrated competitive disruption. This imposes an additional \emph{type-2 thinning} rate $r > 0$ exclusively on the instantaneous leader $\mu_{\pi^{(r)}(1,t)}$, moving it from active status $A$ to restricted status $B$ (from which it still suffers background attrition to $C$ at rate $1$).
  \item \textbf{Network-effect incumbency ($-1 < r < 0$):} Conversely, in markets where strong network effects or switching costs insulate the \#1 platform from churn, the leader enjoys a reduced effective exit rate $1 + r \in (0,1)$ compared to the unit background attrition rate of subdominant rivals.
\end{itemize}
For $r \geq 0$, the transition rates among the states $A$, $B$, and $C$ are summarized by:
\[
  B \xleftarrow[\text{leader}]{r} A \xrightarrow{\quad 1 \quad} C, \qquad B \xrightarrow{\quad 1 \quad} C.
\]
We return to the formal mathematical specification in Section~\ref{subsec:math-description}. As we show below, the stationary market-share distributions under this dynamics are precisely the ``effective fractional balayage'' laws $\mathrm{PD}^{(r)}(m,0)$ for all $r >-1$.

{\em Henceforth, we will abbreviate the phrase ``2-type thinning with alternative leader rate $1+r$''
as ``2TTwALR$(1+r)$.''}

\subsection{Mathematical description of 2TTwALR(1+r)}
\label{subsec:math-description}
Let $\mathcal{F}^{\nearrow}(\N,\N)$ denote the set of increasing functions from $\N$ to $\N$.
This is a monoid (an algebraic semigroup under composition, with identity).
It can be equipped with the projective topology with respect to point evaluations, which makes it 
Polish. Then there is the notion of a continuous convolution semigroup.
If we have a $\mathcal{F}^{\nearrow}(\N,\N)$-valued stochastic process, with sample path $\pi(\cdot,t) : \N \to \N$
defined for $t\geq 0$, let us say that this represents a continuous convolution semigroup if
\begin{itemize}
    \item $\lim_{t \to 0^+} \pi(n,t)=n$, almost surely, for each $n \in \N$, and
    \item we have $\mathbb{P}(\pi(1,t+s)=k_1\, ,\ \dots\, ,\ \pi(n,t+s)=k_n)$ equals
    $$
    \sum_{j_1,\dots,j_n \in \N} \mathbb{P}(\pi(1,t)=j_1\, ,\ \dots\, ,\ \pi(n,t)=j_n)
    \mathbb{P}(\pi(j_1,s)=k_1\, ,\ \dots\, ,\ \pi(j_n,s)=k_n)
    $$
    for every $n \in \N$ and every $k_1,\dots,k_n \in \N$.
\end{itemize}
Of course, we assign probability 0 unless $k_1<\dots<k_n$.
And the summation could be restricted to those $j_1<\dots<j_n$ such that $j_1\leq k_1$, \dots, $j_n\leq k_n$.

Given such a process satisfying the above condition, let us denote by $\mathcal{L}$ the infinitesimal
transition rates (at the identity) or generating functional, defined by
$$
\mathcal{L}_n(k_1,\dots,k_n)\, =\, \lim_{t \to 0^+} \frac{\mathbb{P}(\pi(1,t)=k_1,\dots,\pi(n,t)=k_n)-\delta_{1,k_1} \cdots \delta_{n,k_n}}{t}\, ,
$$
for all $n \in \N$ and all $k_1,\dots,k_n \in \N$ (satisfying $k_1<\dots<k_n$).
For the 2TTwALR$(1+r)$, we denote the special choice
of generating functional as $\Lambda^{(r)}$, and define it as  
$$
\Lambda^{(r)}_n(k_1,\dots,k_n)\, =\, 
\begin{cases}
- (n+r) & \text{ if $k_j=j$, for each $j=1,\dots,n$}\\
1 & \text{ if  $k_j=j+\mathbbm{1}_{\{\ell,\dots,n\}}(j)$ for some $\ell \in \{2,\dots,n\}$, and}\\
1+r & \text{ if $k_j=j+1$ for each $j=1,\dots,n$.}
\end{cases}
$$
We do not prove non-explosion, here. But we believe it can be proved using the method as in Chapter 1, Theorem 3.9 of Liggett \cite{Liggett}.

\begin{lemma}
\label{lem:injection}
    For the continuous convolution semigroup determined by the infinitesimal transition rates $\Lambda^{(r)}$,
    if we denote the sample path as $\pi^{(r)}(\cdot,t) : \N \to \N$,
    defined for all $t\geq 0$, then 
    $$
\mathbb{P}\Big( 
\pi^{(r)}(1,t)\, =\, k_1\, ,\ \pi^{(r)}(2,t)\, =\, k_2\, ,\ \dots\, ,\ 
\pi^{(r)}(N,t)\, =\, k_N\Big)\, =\, \widetilde{Q}^{(r+1)}(k_1,t) \prod_{\ell=2}^N Q(k_\ell-k_{\ell-1},t)\, ,
$$
for
$$
Q(k,t)\, =\, e^{-t} (1-e^{-t})^{k-1}\ \text{ and }\
\widetilde{Q}^{(r+1)}(k,t)\, =\, \frac{\Gamma(r+k)}{\Gamma(k)\Gamma(r+1)}\, e^{-(r+1)t}(1-e^{-t})^{k-1}\, .
$$
\end{lemma}
We will call $\pi^{(r)}(\cdot,t)$ the ``random thinning injection,''
for the 2TTwALR$(1+r)$.

Obviously with this definition $Q=\widetilde{Q}^{(1)}$, and this is the geometric distribution (on $\N$).
Let us define $Q^0(k,t)=\delta_{k,0}$ for all $t\geq 0$.
Then let us also define $Q^R(\cdot,t)$ for each $R \in \N$ by the usual matrix-multiplication
$$
Q^{R}(k_{R},t)\, =\, \sum_{k_1=1}^{k_{R}-R+1}\, 
\sum_{k_2=k_1+1}^{k_{R}-R+2}\,
\cdots\, 
\sum_{k_{R-1}=k_{R-2}+1}^{k_{R}-1} Q(k_1,t)Q(k_2-k_1,t)\cdots Q(k_{R}-k_{R-1},t)\, .
$$
Then by the negative binomial distribution (with the U.~S.~convention)
we see that $Q^{R+1}(\cdot+R,t) = \widetilde{Q}^{(R+1)}(\cdot,t)$.
(This is closely related to whole-number balayage.)
Then the Chu-Vandermonde identity implies that, for general $r,\rho \geq 0$, we have
\begin{equation}
\label{eq:ChuVandermonde}
\sum_{k=1}^{\ell} \widetilde{Q}^{(r+1)}(k,t) \widetilde{Q}^{(\rho+1)}(\ell+1-k,t)\,
=\, \widetilde{Q}^{(r+\rho+2)}(\ell,t)\, ,
\end{equation}
for each $\ell \in \N$.
See, for instance, Andrews, Askey and Roy, Corollary 2.2.3 \cite{AndrewsAskeyRoy}.

Now let us finish the proof of the lemma.

\begin{proof}
We begin by showing that
\begin{equation}
\label{eq:MarkovFirst}
\sum_{N=1}^{\ell} \mathbb{P}\Big(\pi^{(r)}(N,t)\, =\, \ell\Big)
\mathbb{P}\Big(\pi^{(r)}(1,s)\, =\, N\Big)\, =\, 
\mathbb{P}\Big(\pi^{(r)}(1,s+t)\, =\, \ell\Big)\, .    
\end{equation}    
By (\ref{eq:ChuVandermonde}), the desired identity is equivalent to
$$
\sum_{k} \widetilde{Q}^{(r+1)}(k,t) \widetilde{Q}^{(\rho+1)}(\ell+1-k,t)\,
=\, \widetilde{Q}^{(r+\rho+2)}(\ell,t)\, .
$$
That is what we must show to establish (\ref{eq:MarkovFirst}).
But note that the left-hand-side is actually
$$
\sum_N \frac{\Gamma(r+N-1+\ell-N+1)}{\Gamma(\ell-N+1)\Gamma(r+N)}\,
q^{r+N}(1-q)^{\ell-N}
\frac{\Gamma(r+N)}{\Gamma(N)\Gamma(r+1)}\, 
z^{r+1} (1-z)^{N-1}\, ,
$$
for $q=e^{-t}$ and $z=e^{-s}$.
This simplifies to
$$
\frac{\Gamma(r+\ell)}{\Gamma(\ell)\Gamma(r+1)} 
(qz)^{r+1}
\sum_N \frac{(\ell-1)!}{(\ell-N)!(N-1)!}\, 
(1-q)^{\ell-N}\big(q(1-z)\big)^{N-1}\, .
$$
And by the binomial formula, this equals
$$
\frac{\Gamma(r+\ell)}{\Gamma(\ell)\Gamma(r+1)}\, (qz)^{r+1} (1-qz)^{\ell-1}\, ,
$$
which is $\widetilde{Q}^{(r+1)}(\ell,s+t)$.
So (\ref{eq:MarkovFirst}) has been established.

Now, by its definition, the gaps for the 2TTwALR$(1+r)$ are independent
of the first value and are distributed as IID Geometric$(e^{-t})$ random variables.
This is reflected in the product of the $Q(\cdot,t)$-kernels.
\end{proof}
It is an exercise left for the reader that with this definition we do have the infinitesimal transition rates $\mathcal{L}$ being equal to $\Lambda^{(r)}$.

\begin{definition}
    If we have a dust-free random partition structure with infinitely many points $(\mathsf{X}_1,\mathsf{X}_2,\dots)$, ordered
    so that $\mathsf{X}_1>\mathsf{X}_2>\dots>0$, then the 2TTwALR$(1+r)$ is defined by letting $\pi^{(r)}_t$
    be the associated random thinning injection,
    sampled independently of $(\mathsf{X}_1,\mathsf{X}_2,\dots)$. Then let 
    $$
    \mathsf{Y}^{(r)}(n,t)\, =\, \mathsf{X}_{\pi^{(r)}(n,t)}\, ,
    $$
    for each $n \in \N$. And let $\mathsf{C}^{(r)}(t)=-\ln(\mathsf{Y}^{(r)}(1,t)+\mathsf{Y}^{(r)}(2,t)+\dots)$.
    Finally define
    $$
    \widehat{\mathsf{Y}}^{(r)}(n,t) = \exp\left(\mathsf{C}^{(r)}(t)\right) \mathsf{Y}^{(r)}(n,t)\, ,
    $$
    for each $n \in \N$.
    Then $(\widehat{\mathsf{Y}}^{(r)}(1,t),\widehat{\mathsf{Y}}^{(r)}(2,t),\dots)$ is the 2TTwALR$(1+r)$.
\end{definition}

\begin{remark}
    In Section \ref{sec:Graphical} we will establish a graphical representation of this thinning along the lines of Harris \cite{Harris}.
\end{remark}

In the next  section, we will characterize those random partition structures that are invariant
with respect to the 2TTwALR$(1+r)$.
For now, let us observe the following corollary of Lemma \ref{lem:injection}.

\begin{cor}
\label{cor:Gamma}
    Suppose that $\xi_1<\xi_2<\dots$ is a homogeneous Poisson point process with intensity 1 on $(0,\infty)$.
    Then for $r\geq 0$ let $(\widetilde{\xi}^{(r)}(1,t),\widetilde{\xi}^{(r)}(2,t),\dots)$
    be defined as 
    $$
        \widetilde{\xi}^{(r)}(n,t)\, =\, \xi_{\pi^{(r)}(n,t)}\, ,
    $$
    where $\pi^{(r)}(\cdot,t)$ is the random thinning injection for 2TTwALR$(1+r)$, sampled
    independently of $(\xi_1,\xi_2,\dots)$.
    In this case, 
    $$
    \mathbb{E}\left[\exp\left(-\lambda \widetilde{\xi}^{(r)}(1,t)\right)\right]\, 
    =\, \frac{(1+\lambda)^r}{(1+\lambda e^{t})^{1+r}}\, .
    $$
    Of course, $\widetilde{\xi}^{(r)}(n,t)-\widetilde{\xi}^{(r)}(n-1,t)$
    for $n\geq 2$ are all independent of $\widetilde{\xi}^{(r)}(1,t)$
    and are IID Geometric($e^{-t}$) random variables.
\end{cor}
In particular, the Laplace transform of $e^{-t} \widetilde{\xi}^{(r)}(1,t)$
converges to $(1+\lambda)^{-(1+r)}$ so that the random variable
$e^{-t} \widetilde{\xi}^{(r)}(1,t)$ converges in-distribution
to a Gamma-$(1+r,1)$ random variable.

\begin{proof}
    This follows because the probability generating function for $\widetilde{Q}^{(r+1)}(\cdot,t)$ is
    $$
    \sum_{k=1}^{\infty} \widetilde{Q}^{(r+1)}(k,t) z^k\, =\, z \left(\frac{q}{1-(1-q)z}\right)^{r+1}\, ,
    $$
    for $q=e^{-t}$. Then the result follows from the chain-rule for Laplace transforms
    of functions
    $$
    \mathbb{E}\left[\exp\left(-\lambda \widetilde{\xi}^{(r)}(1,t)\right)\right]\, 
    =\, \sum_{k=1}^{\infty} \widetilde{Q}^{(r+1)}(k,t) \left(\frac{1}{1+\lambda}\right)^k\, =\, \frac{1}{1+\lambda}\, \left(\frac{1+\lambda}{1+q^{-1}\lambda}\right)^{r+1}\, .
    $$
\end{proof}
Because of all this, $e^{-t}\widetilde{\xi}^{(r)}(1,t)$ converges in distribution
to the sum of an exponential-1 and an independent Gamma-$(r,1)$,
while all the increments $e^{-t}\widetilde{\xi}^{(r)}(n,t)-e^{-t}\widetilde{\xi}^{(r)}(n-1,t)$,
for $n=2,3,\dots$ are independent of $e^{-t}\widetilde{\xi}^{(r)}(1,t)$,
and they converge in distribution to IID exponential-1 random variables.
If we post-compose with $\Phi_q$, for $m=1/q$, and we normalize, then we obtain the set-up of $\mathrm{PD}^{(r)}(m,0)$.

\section{Invariance for 2TTwALR$(r)$}
An easy fact is the following.
\begin{lemma}
    If $(\mathsf{X}_1,\mathsf{X}_2,\dots)$ is distributed as a $\mathrm{PD}(m,0)$ for $0<m<1$,
    then the 1-type thinning -- that is the 2TTwALR$(1)$ --
    defined as $(\widehat{\mathsf{Y}}^{(0)}_1(t),\widehat{\mathsf{Y}}^{(0)}_2(t),\dots)$ 
    is $t$-stationary. Moreover, $\lim_{t \to \infty} t^{-1} \mathsf{C}^{(0)}(t) = 1/m$, a.s.
\end{lemma}
This fact is well-known. It is a special case of a more general fact for
$\mathrm{PD}(m,0)$ partition structures.
For the more general quasi-stationary property, see Ruelle's article \cite{Ruelle}.
We also wish to cite Aizenman and Ruzmaikina \cite{AizenmanRuzmaikina}.


Now note that $\pi^{(r+1)}_t$ is distributed identically to $\pi^{(r)}_t(\cdot+1)-1$.
In an infinitesimal time interval if either the first or second particle
is thinned, then the second particle moves.
This is 1-point balayage.
For any $R \in \{0,1,\dots\}$ we can deduce that $\mathrm{PD}^{(R)}(m,0)$
is a random partition structure which is invariant for 2TTwALR$(1+R)$.



\begin{prop}
    Suppose $r>-1$. For any $q>1$, let $(\mathsf{X}_1,\mathsf{X}_2,\dots)$ be a 
    $\mathrm{PD}^{(r)}(1/q,0)$ random
    partition structure. Let $(\widehat{\mathsf{Y}}^{(r)}_1(t),\widehat{\mathsf{Y}}^{(r)}_2(t),\dots)$ be the 2TTwALR$(1+r)$ derived from $(\mathsf{X}_1,\mathsf{X}_2,\dots)$.
    Then we have time stationarity of 
    $(\widehat{\mathsf{Y}}^{(r)}_1(t),\widehat{\mathsf{Y}}^{(r)}_2(t),\dots)$.
    Moreover, $\lim_{t \to \infty} t^{-1} \mathsf{C}^{(r)}(t)=q$.
\end{prop}

We will establish this in Section \ref{sec:Graphical} by using a graphical representation, which we believe may also be of interest.
In the next two subsections, we characterize those random partition structures which are invariant
with respect to 2TTwALR$(r)$ (among dust-free random partition structures with infinitely many points).

\subsection{Perturbation of stationary random partition structures}

Let $\mathcal{M}_{+}^{(1)}(\R)$ denote the set of all measures $\mu$ on $\R$
such that 
\begin{equation}
\label{eq:condition}
\lim_{R \to \infty} \frac{\mu(-R,R)}{2R}\, =\, 1\, .
\end{equation}
Now let $\Delta$ denote the set of probability measures $\mathfrak{P}$ on $\mathcal{M}_{+}^{(1)}(\R)$
such that if $M$ is a sample from $\mathfrak{P}$ then 
$M(\cdot+a)$ is equal-in-distribution to $M$ for each $a \in \R$.
So the distribution $\mathfrak{P}$ is translation invariant as a measure.
For example, consider $M_{\omega}(dx)=2\cos^2(x+\omega)\, dx$ for $\omega$ uniformly distributed on 
$\R/(2\pi \Z)$. Then the resulting distribution $\mathfrak{P}$ is translation invariant.

Then $\Delta$ forms a Poulsen simplex.
See for instance \cite{LOS}.
(Also see \cite{Kallenberg}.)
Given an $M\in \mathcal{M}_{+}^{(1)}(\R)$ let us define $G_M : \R \to (0,\infty)$ by
$$
G_M(x)\, =\, \int_{-\infty}^x e^y\, M(dy)\, .
$$
Note that defining ${M}'$ by ${M}'(dy) = e^y\, M(dy)$,
we have that ${M}'(\cdot+a)$ is equal in distribution to $e^a {M}'$ for each $a \in \R$.
Hence, $G_M(\cdot+a)$ is equal in distribution to $e^a G_M$.
Finally, consider the mapping $\mathcal{H}_M : (0,\infty) \to (0,\infty)$
by $\mathcal{H}_M(t) = G_M(\ln(t))$.
Then $\mathcal{H}_M(e^a\cdot)$ is equal in distribution to $e^a \mathcal{H}_M$.
For example,
$$
\mathcal{H}_{M_{\omega}}(t)\, =\, t \left[1-a + 2a \sin^2(\tau+\omega+\ln(t))\right]\, ,
$$
from the example before, where $0<a<1$ and $\tau \in [0,2\pi)$ is a phase shift.

\begin{definition}
    If we have a dust-free random partition structure $(\mathsf{X}_1,\mathsf{X}_2,\dots)$
    with infinitely many points
    and if we have a measure $\mathfrak{P} \in \Delta$, let $M \in \mathcal{M}_+^{(1)}(\mathbb{R})$ 
    be a random sample from $\mathfrak{P}$, independent of $(\mathsf{X}_1,\mathsf{X}_2,\dots)$.
    Then we say that the $\mathfrak{P}$-perturbed random partition structure is the distribution of
    $(\sigma_M(\mathsf{X}_{\cdot})_1,\sigma_M(\mathsf{X}_{\cdot})_2,\dots)$ where
    $$
    \sigma_M(\mathsf{X}_{\cdot})_n\, =\, \frac{\mathcal{H}_M(\mathsf{X}_n)}{\sum_{\nu=1}^{\infty} \mathcal{H}_M(\mathsf{X}_\nu)}\, ,
    $$
    for all $n \in \N$. 
\end{definition}

If $\textsc{Law}$ is the law of $(\mathsf{X}_1,\mathsf{X}_2,\dots)$, then let us write
$\textsc{Pert}_{\mathfrak{P}}(\textsc{Law})$ for the law of 
$(\sigma_M(\mathsf{X}_{\cdot})_1,\sigma_M(\mathsf{X}_{\cdot})_2,\dots)$.

\begin{lemma}
    Suppose that $\textsc{Law}$ is the law for the dust-free random partition structure
    $(\mathsf{X}_1,\mathsf{X}_2,\dots)$ with infinitely many parts.
    If $\textsc{Law}$ is invariant for 2TTwALR$(1+r)$
then for any $\mathfrak{P} \in \Delta$, we have that
$\textsc{Pert}_{\mathfrak{P}}(\textsc{Law})$
    is also invariant for 2TTwALR$(1+r)$.
\end{lemma}

\begin{proof}
    By invariance, $(\widehat{\mathsf{Y}}^{(r)}(1,t),\widehat{\mathsf{Y}}^{(r)}(2,t),\dots,)$ is equal in distribution
    to $(\mathsf{X}_1,\mathsf{X}_2,\dots)$.
    Hence, the rescaled process $(e^{\mathsf{C}^{(r)}(t)}\mathsf{Y}^{(r)}(1,t),e^{\mathsf{C}^{(r)}(t)}\mathsf{Y}^{(r)}(2,t),\dots)$ is equal in distribution 
    to $(\mathsf{X}_1,\mathsf{X}_2,\dots)$.
    Then by the properties of $\mathfrak{P}$, that means that $(\mathcal{H}_M(\mathsf{Y}^{(r)}(1,t)),
    \mathcal{H}_M(\mathsf{Y}^{(r)}(1,t)),\dots)$ is equal to the scaled sequence
    $$
    (e^{-\mathsf{C}^{(r)}(t)}\mathcal{H}_{\widetilde{M}}(\widehat{\mathsf{Y}}^{(r)}(1,t)),
    e^{-\mathsf{C}^{(r)}(t)}\mathcal{H}_{\widetilde{M}}(\widehat{\mathsf{Y}}^{(r)}(2,t)),\dots)\, ,
    $$
    where $\widetilde{M}$ is equal in distribution to $M$.
    But when we normalize, we have that 
    $$
    \frac{\sigma_M(\mathsf{X}_{\cdot})_{\pi^{(r)}(n,t)}}
    {\sum_{\nu=1}^{\infty}\sigma_M(\mathsf{X}_{\cdot})_{\pi^{(r)}(\nu,t)}}\, =\, \frac{\mathcal{H}_M(\mathsf{X}_{\pi^{(r)}(n,t)})}{\sum_{\nu=1}^{\infty} \mathcal{H}_M(\mathsf{X}_{\pi^{(r)}(\nu,t)})}\, ,
    $$
    by cancellation of the normalizers.
    And this is equal to 
    $$
    \frac{\mathcal{H}_M(\mathsf{Y}^{(r)}(n,t))}{\sum_{\nu=1}^{\infty} \mathcal{H}_M(\mathsf{Y}^{(r)}(\nu,t))}\, .
    $$
    And this is equal in distribution to $\mathcal{H}_M(\widehat{\mathsf{Y}}^{(r)}(n,t))/\sum_{\nu=1}^{\infty}
    \mathcal{H}_M(\widehat{\mathsf{Y}}^{(r)}(\nu,t))$, as a $n$-indexed sequence, just because of cancellation
    of the multiplier $e^{-\mathsf{C}^{(r)}(t)}$ from the numerator and denominator.
    But then since $\widehat{\mathsf{Y}}^{(r)}(\cdot,t)$ is equal in distribution to
    $\mathsf{X}_{\cdot}$ (by the invariance property with respect to 2-type thinning
    with leader rate $r$), we see that the result
    is equal in distribution to $\mathcal{H}_M(\mathsf{X}_n)/\sum_{\nu=1}^{\infty} \mathcal{H}_M(\mathsf{X}_{\nu})$, as an $n$-indexed sequence.
    But since the latter is $\sigma_M(\mathsf{X}_{\cdot})$, we have reached the desired result.
\end{proof}
This elementary result would also hold if $\mathfrak{P}$ had the property that $\mathcal{H}_M(e^a\cdot)$
were equal in distribution to $ e^{t a}\mathcal{H}_M$ for all $a$ and some (fixed) $t>0$.
But we choose $t=1$, which is encoded in the condition (\ref{eq:condition}), to remove
the extra parameter of non-uniqueness.

\subsection{Consequence of the ergodic theorem}

The following result follows from Birkhoff's ergodic theorem.
\begin{prop}
\label{prop:ergodic}
    Consider a dust-free 
    random partition structure $(\mathsf{X}_1,\mathsf{X}_2,\dots)$
    with infinitely many points, and  let 
    $$
    (\widehat{\mathsf{Y}}^{(r)}(1,t),\widehat{\mathsf{Y}}^{(r)}(2,t),\dots)
    $$
    be the 2-type thinning with additional leader rate $r\geq 0$.
Suppose that $(\widehat{\mathsf{Y}}^{(r)}(1,t),\widehat{\mathsf{Y}}^{(r)}(2,t),\dots)$
    is $t$-stationary.
Suppose that $\mathbb{E}[\mathsf{C}^{(r)}(t)]<\infty$. Then
$$
\lim_{t \to \infty} \frac{\mathsf{C}^{(r)}(t)}{t}\, =\, \mathsf{Q}\, ,
$$
almost surely for a  random variable $\mathsf{Q}$, which is measurable with respect to the shift-invariant \(\sigma \)-algebra
\end{prop}

For the proof of Proposition \ref{prop:ergodic}, we refer to any standard theorem of the ergodic theorem.
For example, Billingsley's Theorem 24.1 \cite{Billingsley} or Dudley's Theorem 8.4.1 \cite{Dudley}.

By the nature of thinning, we must have $\mathsf{Q}>0$, almost surely.

\begin{theorem}
\label{thm:characterization}
    Consider a dust-free random partition structure $(\mathsf{X}_1,\mathsf{X}_2,\dots)$ with infinitely
    many points, and  let 
    $$
    (\widehat{\mathsf{Y}}^{(r)}(1,t),\widehat{\mathsf{Y}}^{(r)}(2,t),\dots)
    $$
    be the 2TTwALR$(1+r)$ for $r>-1$.
Suppose that $(\widehat{\mathsf{Y}}^{(r)}(1,t),\widehat{\mathsf{Y}}^{(r)}(2,t),\dots)$
    is $t$-stationary, and assume $\mathbb{E}[\mathsf{C}^{(r)}(t)]<\infty$.
    Then $\mathsf{Q}=\lim_{t \to \infty} \mathsf{C}^{(r)}(t)/t$ is strictly greater than 1, a.s.
    Moreover, conditioning on $\mathsf{Q}=q$ (by using a regular conditional distribution)
    we have that $(\mathsf{X}_1,\mathsf{X}_2,\dots)$ is 
    $\textsc{Pert}_{\mathfrak{P}(q)}(\mathrm{PD}^{(r)}(1/q,0))$ for a unique $\mathfrak{P}(q) \in \Delta$.
\end{theorem}

The proof of this theorem is  similar to the proof of an analogous result from
an article of Vermesi, Wei and one of the authors \cite{SVW}.
The main difference is that the rescaled limit of $\pi^{(r)}(\cdot,t)$ converges to the 
distribution we have been describing: $\xi_1$ is Gamma-$(1+r,1)$, and conditional on $\xi_1$
the points $\xi_2,\xi_3,\dots$ are a homogeneous Poisson point process of rate 1 on $(\xi_1,\infty)$.
Having said this, the details are so similar that we omit them.

\section{The graphical representation}
\label{sec:Graphical}

\begin{figure}
    \centering
    \includegraphics[width=0.75\linewidth]{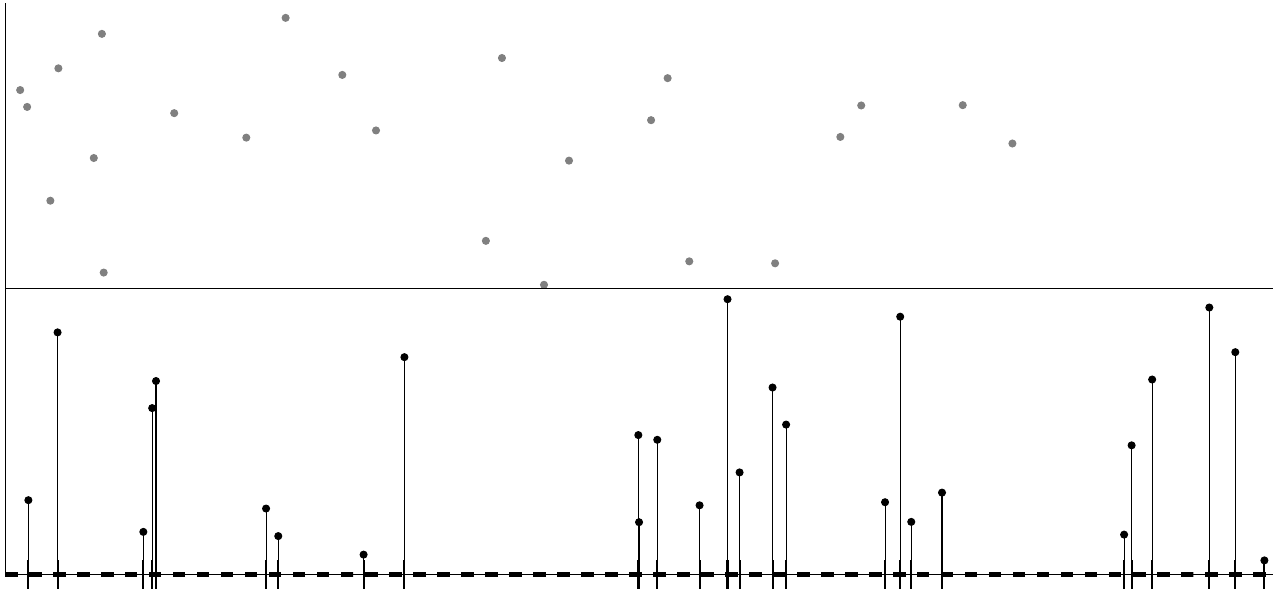}
    \caption{The graphical representation for 1-type thinning: All the dots represent the intensity 1 two-dimensional
    Poisson point process $\Xi$ on the positive quadrant. 
    The vertical hatch marks on the horizontal axis represent $\Pi(\Xi \cap \mathcal{R})$.
    }
    \label{fig:placeholder0}
\end{figure}

We may construct a  graphical representation.
We begin with the graphical representation for 1-type thinning.

Let $\Xi$ be a standard
Poisson process on $(0,\infty)\times (0,\infty)$,
with 2d (area) Lebesgue measure as the intensity,
viewed as  a random subset $\Xi \subset (0,\infty)\times (0,\infty)$.
Then considering the mapping
$$
\mathcal{J}_t(x,y)\, =\, (e^{-t}x,e^ty)\, ,
$$
which is a linear hyperbolic evolution, we know that the law of $\mathcal{J}_t(\Xi)$ is equal to the law of $\Xi$ because 
$\mathcal{J}_t$ preserves 2d Lebesgue measure.

So now, let $\Pi : \R^2 \to \R$ just be the projection $\Pi(x,y) = x$.
Also, let $\mathcal{R}$ be the semi-infinite rectangle $\mathcal{R} = \{(x,y)\, :\, x \in [0,\infty)\, ,\ y\in(0,1)\}$.
Then define $(\xi_1(t),\xi_2(t),\dots)$ such that $\xi_1(t)<\xi_2(t)<\dots$
and
$$
\{\xi_n(t)\, :\, n \in \N\}\, =\, \Pi(\mathcal{J}_t(\Xi)\cap \mathcal{R})\, .
$$
Then it is easy to see that $(\xi_1(t),\xi_2(t),\dots)$ is identical as a process to the rate-1 thinning of $(\xi_1(0),\xi_2(0),\dots)$ followed
by rescaling by $e^{-t}$ to keep the intensity fixed.
The thinning is accomplished by taking $y$ to $e^ty$ which is less 
than 1 if and only if $y<e^{-t}$. So this is thinning of points in $\mathcal{R}$ at rate 1.
See Figure \ref{fig:placeholder0} for an example of the Poisson process $\Xi$ and $\Pi(\Xi\cap \mathcal{R})$.

This gives a graphical representation for uniform thinning.
This is a well-known fact. A good reference is Aldous \cite{Aldous2003}.
In a paper by Vermesi, Wei and one of the current authors \cite{SVW}
this precise representation was described.

This accomplishes the standard 1-type thinning, at the level of a Poisson
point process.
I.e., it implements 2TTwALR$(1)$. One can then apply $\Theta_q$ to show that the Poisson-Dirichlet $\mathrm{PD}(1/q,0)$-model is thinning-invariant for the 1-type
thinning (IID thinning of all parts, with no additional rate for the leader, $r=0$).

\subsection{2TTwALR(1+r) for positive r by boundary ``sinks"}

\begin{figure}
    \centering
    \includegraphics[width=\linewidth]{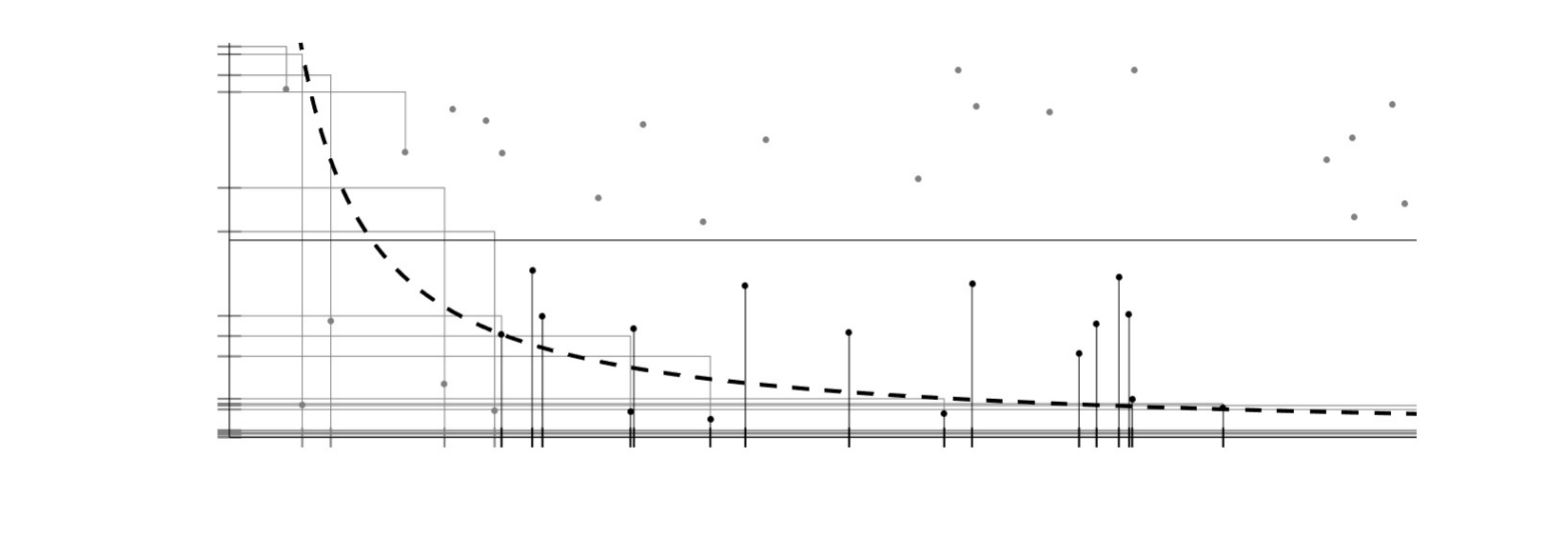}
    \caption{The graphical representation for $2TTwALR(1+r)$ for $r=3$: All the dots represent the intensity 1 two-dimensional
    Poisson point process $\Xi$ on the positive quadrant. The horizontal hatch marks on the vertical
    axis represent the points of $\Theta^{(r)}$. The horizontal line at height 1 represents the boundary
    of $\mathcal{R}$. The gray elbow connections represent the matchings of points in $\Theta^{(r)}$
    to points in $\Omega^{(r)} \subset \Xi$. For points in $\Xi \cap \mathcal{R}$, we draw verticals down towards vertical hatch
    marks on the horizontal axis. The $\Pi$-images of $\Omega^{(r)} \cap \mathcal{R}$ are indicated by gray hatch marks because they are {\em not} part of $(\xi_1^{(r)}(t),\xi^{(r)}_2(t),\dots)$. Those are the points which are effectively ``balayaged'' 
    in the fractional balayage.
    The vertical hatch marks on the horizontal axis which are black are the points of the resulting point process that survives after
    the effective fractional balayage. The black dashed line is the ``fluid limit'' approximation
    of the boundary of the matched points $xy=r$.}
    \label{fig:placeholder}
\end{figure}

\begin{figure}
    \centering
    \includegraphics[width=\linewidth]{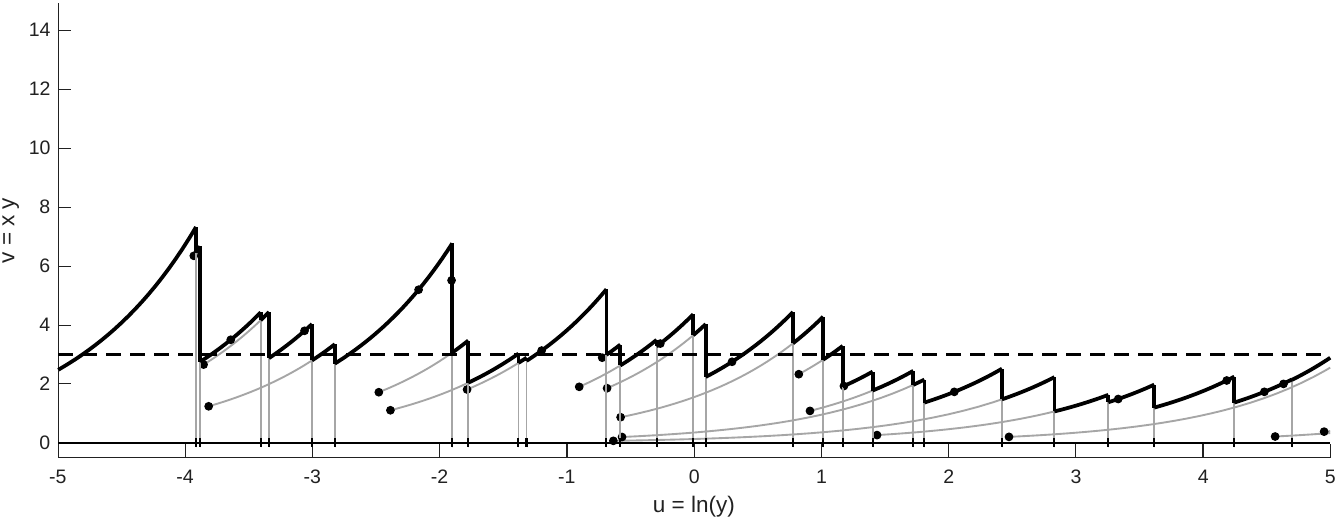}
    \caption{The graphical representation in the  coordinates $u=\ln(y)$ and $v=xy$:
    This time we only plot the matched points. The black dashed line is the ``fluid limit''
    for the boundary of this region $v=r$. The solid black curve is the boundary of the elbow
    connectors implementing the matching.}
    \label{fig:placeholder2}
\end{figure}

Now we augment our graphical representation.

Let $r\,dy/y$ be a measure on the $y$-axis ($x=0$), which is scale invariant.
For technical reasons, let us define
$$
\mu^{(T)}(dy)\, =\, \mathbf{1}_{(0,e^T)}(y)\, \frac{dy}{y}\, ,
$$
for each $T \in (0,\infty)$, meaning we started our process in the past at time $-T$.
But we are really interested in $\mu^{(\infty)}$, which we obtain by taking $T \to \infty$.
Even though $x=0$, the measure $\delta_0 \otimes \mu^{(\infty)}$ is invariant under the action of $\mathcal{J}_t$, since the scale-invariant measure is invariant under the rescaling $y \mapsto e^t y$.
Of course, $(\delta_0\otimes \mu^{(T)}) \circ \mathcal{J}_t^{-1} = \delta_0 \otimes \mu^{(T+t)}$.

Let $\Theta^{(r,T)}$ be a random subset of $(0,\infty)$ corresponding to a Poisson point process on $(0,\infty)$ with intensity $r\cdot \mu^{(T)}$. Then, we may enumerate the points of $\Theta^{(r,T)}$ as 
$(\theta^{(r,T)}_1,\theta^{(r,T)}_2,\dots)$ such that
\[
e^T>\theta_{1}^{(r,T)} > \theta_2^{(r,T)} > \theta_3^{(r,T)} > \cdots
\]
such that $\{\theta_n^{(r)} : n \in \mathbb{Z}\} = \Theta^{(r)}$. Then proceeding inductively along the 
indices, we match each $\theta_n^{(r,T)}$ to the  point $(x_n,y_n) \in \Xi$ 
satisfying
$$
(x_n,y_n)\, \in \, \mathcal{V}_n^{(r,T)}\, \stackrel{\mathrm{def}}{:=}\,
\{(x,y) \in \Xi \setminus \mathcal{W}^{(r,T)}_{n-1}\, :\, y<\theta_n^{(r,T)}\}\, ,
$$
and $x_n = \min(\Pi(\mathcal{V}_n^{(r,T)}))$, where $\mathcal{W}_{0}^{(r,T)}=\emptyset$
and $\mathcal{W}_n^{(r,T)} = \{(x_1,y_1),\dots,(x_n,y_n)\}$ for each $n \in \N$.
For finite $T$, this greedy matching is almost surely well-defined and non-degenerate: 
$\mathcal{W}_{\infty}^{(r,T)} := \bigcup_{n=1}^{\infty} \Omega^{(r,T)}_n$ is not all of $\Xi$.

Even as $T \to \infty$, it is true that the matching almost surely does not exhaust $\Xi$. This follows from Corollary \ref{cor:Gamma}.
Define $\Omega^{(r,T)}$ to be the set
$$
\Omega^{(r,T)}\, =\, \{(x_n,y_n)\, :\, n \in \N\, ,\ \theta_n^{(r,T)} > 1\, ,\ y_n<1\}\, .
$$
Then the left-most point of $(\Xi \cap \mathcal{R})\setminus \Omega^{(r,T)}$
is at position $\widetilde{\xi}^{(r)}(n,T)$ from Corollary \ref{cor:Gamma}, which is stochastically dominated by a 
Gamma-$(1+r,1)$ random variable. (Hence the distributions are tight even as $T \to \infty$.)
So we may take the $T \to \infty$ limit.

An interesting question is to study $\mathfrak{M}^{(r,T)}$, the cardinality of $\Omega^{(r,T)}$.
Let us define $\mathfrak{M}^{(r,T)}_t$ to be the analog of $\mathfrak{M}^{(r,T+t)}$
when we evolve $[0,\infty)\times (0,\infty)$ according to $\mathcal{J}_t$.
It is easy to see that this is a birth-and-death process with easily calculable rates, as in the next lemma.
\begin{lemma}\label{lem:mminf}
For the birth-and-death process on $\{0, 1, \dots\}$ (equivalently, the customer-count process of an $M/M/\infty$ queue with arrival rate $r$ and unit service rate) with generator rates
\[
\lim_{\epsilon \to 0^+} \frac{\mathbb{P}(\mathfrak{M}_{t+\epsilon} = k\, |\, \mathfrak{M}_t = j) - \delta_{j,k}}{\epsilon} =
\begin{cases}
-(j + r) & \text{if } k = j, \\
j & \text{if } k = j - 1, \\
r & \text{if } k = j + 1,
\end{cases}
\]
the unique invariant measure is $\mathrm{Poisson}(r)$.
\end{lemma}

\begin{proof}
The detailed balance equations $r\,\pi(j) = (j+1)\,\pi(j+1)$ for $j \geq 0$ uniquely determine $\pi(j) = \frac{r^j}{j!}\pi(0)$, and normalizing $\sum_{j=0}^\infty \pi(j) = 1$ yields $\pi(j) = e^{-r} r^j / j!$. Because the chain is irreducible and non-explosive on $\{0,1,\dots\}$, $\mathrm{Poisson}(r)$ is its unique stationary distribution.
\end{proof}

Note that $\mathfrak{M}^{(r,T)}$ is identical in distribution to $\mathfrak{M}^{(r,0)}_T$.
Since $\mathcal{M}^{(r,0)}_0=0$, a.s., we see that $\mathfrak{M}^{(r,T)}$ is stochastically dominated
by a Poisson-$r$ random variable, for all $T \in (0,\infty)$.
This is another way to see that the $T \to \infty$ limit of the graphical representation is well-defined.

Henceforth, we will assume $T$ is taken to $\infty$.
We will drop the $\infty$ from the notation: so we write $\Omega^{(r)}$
instead of $\Omega^{(r,\infty)}$.
Let us enumerate $\Theta^{(r)} = \{\theta^{(r)}_n\, :\, n \in \mathbb{Z}\}$
where $\dots >\theta^{(r)}_{-1}>\theta^{(r)}_0>\theta^{(r)}_{1}>\dots$,
and let $\theta^{(r)}_0>1>\theta^{(r)}_1$.
Then, letting $(x_n,y_n) \in \Xi$ be the point matched to $\theta^{(r)}_n$ for each $n \in \mathbb{Z}$.
let us also define 
$$
\Omega^{(r)}_t\, =\, \{(e^{-t} x_n,e^t y_n)\, :\, n \in \mathbb{Z}\, ,\ e^t\theta_n^{(r)}>1\, ,\ 
e^t y_n<1\}\, .
$$
Then we may define $(\xi^{(r)}_1(t),\xi^{(r)}_2(t),\dots)$ such that 
$\xi^{(r)}_1(t)<\xi^{(r)}_2(t)<\dots$
and
$$
\{\xi_1^{(r)}(t)\, ,\ \xi_2^{(r)}(t)\, ,\ \dots\}\, =\, 
\Pi\left(\mathcal{J}_t(\Xi) \cap \mathcal{R}\right) \setminus \Omega^{(r)}_t\, .
$$
By the construction, the process of $(\xi^{(r)}_1(t),\xi^{(r)}_2(t),\dots)$
is stationary in time.
An example is shown in Figure \ref{fig:placeholder}.

Now according to the previous subsection, if any point $(x,y) \in \Xi \cap \mathcal{R}$ is mapped by 
$\mathcal{J}_t$ outside of $\mathcal{R}$, then it is thinned (by the usual
1-type thinning). But additionally, if $(x,y)$ is matched to a point  $\theta_n^{(r)}$, then it is considered thinned as soon as $\mathcal{J}_t$
maps $(0,\theta_n^{(r)})$ outside $\{0\} \times (0,1)$. We claim that this additional
thinning exactly accomplishes the ``effective fractional $r$-balayage'' which
we introduced in the last section by using a Gamma-$(r,1)$
random shift $\mathsf{Z}^{(r)}$. But now we have it at the process-level.

\begin{proposition}\label{prop:sink_gamma}
In the $T \to \infty$ stationary matching above, let $\mathsf{Z}^{(r)}$ denote the rightmost horizontal coordinate in $\Omega^{(r)}$.
Then
$\mathsf{Z}^{(r)}$ has the $\mathrm{Gamma}(r, 1)$ distribution.
\end{proposition}

\begin{proof}
For $r = 0$, $\Theta^{(0)} = \emptyset$ almost surely, so $\mathsf{Z}^{(0)} = 0$ and the claim is immediate. Assume $r > 0$.


By induction along $k = -M, \dots, 0$, each sink $\theta_k^{(r)}$ strictly advances the excavated horizontal frontier from $x_{k-1}$ to $x_k = x_{k-1} + \Delta x_k$, where $\Delta x_k$ is the horizontal distance to the leftmost point of $\Xi$ in the untouched vertical strip $(x_{k-1}, \infty) \times (0, \theta_k^{(r)})$. Since a strip of width $\Delta x$ and height $\theta_k^{(r)}$ has area $\theta_k^{(r)}\,\Delta x$, we have
\[
\mathbb{P}(\Delta x_k > \Delta x \mid \theta_k^{(r)}) = \exp\!\big(-\theta_k^{(r)}\,\Delta x\big),
\]
so $\gamma_k := \theta_k^{(r)}\,\Delta x_k$ are i.i.d.\ $\mathrm{Exponential}(1)$ random variables independent of $\Theta^{(r)}$. Summing over all sinks in $(1, e^T)$ and letting $T \to \infty$ gives the exponential shot-noise (or Vervaat perpetuity) representation
\[
\mathsf{Z}^{(r)} = \sum_{\theta_k^{(r)} > 1} \frac{\gamma_k}{\theta_k^{(r)}} = \sum_{U_k > 0} \gamma_k\,e^{-U_k},
\]
where $\{U_k = \ln \theta_k^{(r)} > 0\}$ is a homogeneous Poisson point process of rate $r$ on $(0, \infty)$ independent of the i.i.d.\ $\mathrm{Exponential}(1)$ marks $\{\xi_k\}$ (cf.\ \cite{GPY2006,Kallenberg}).

By Campbell's exponential formula (for the Laplace functional of Poisson processes) applied to the marked Poisson process $\{(U_k, \gamma_k)\}$ on $(0,\infty) \times (0,\infty)$ with intensity $r\,du \otimes e^{-x}\,dx$, the Laplace transform of $\mathsf{Z}^{(r)}$ for any $\lambda \geq 0$ is
\begin{align*}
\mathbb{E}\!\left[\exp\left(-\lambda \mathsf{Z}^{(r)}\right)\right]
&= \exp\!\left( - \int_0^\infty \int_0^\infty \left(1 - e^{-\lambda x e^{-u}}\right) r\,du\,e^{-x}\,dx \right) \\
&= \exp\!\left( - r \int_0^\infty \left(1 - \frac{1}{1 + \lambda e^{-u}}\right) du \right) \\
&= \exp\!\left( - r \int_0^\infty \frac{\lambda e^{-u}}{1 + \lambda e^{-u}}\,du \right) \\
&= \exp\!\left( - r \Big[ -\ln(1 + \lambda e^{-u}) \Big]_0^\infty \right) = (1 + \lambda)^{-r}.
\end{align*}
Thus $\mathsf{Z}^{(r)} \sim \mathrm{Gamma}(r, 1)$, proving (1).
\end{proof}

Also note that
conditioned on $\mathsf{Z}^{(r)}$, we have
$$
\Omega^{(r)}\, =\, \{(x,y) \in \Xi \cap \mathcal{R}\, :\, x\leq \mathsf{Z}^{(r)}\}\, ,
$$
and $(\Xi \cap \mathcal{R}) \setminus \Omega$ is a homogeneous Poisson point process of intensity 1 on $(\mathsf{Z}^{(R)},\infty) \times (0,1)$.
Consequently, the boundary-sink thinning with $\mathcal{J}_t$ realizes the 2TTwALR$(1+r)$ for
$\Pi(\Xi \cap \mathcal{R})$ at the process level.

The fluid limit is the expectation of $\mathsf{Z}^{(r)}$ mapped forward and backwards according to
$\mathcal{J}_t$. 
(We do not report the details here, but for us a first step was to develop the fluid limit approximation.
In this problem, fortunately, that is closely connected to the real 2TTwALR$(1+r)$ process.)
It is also the expectation for the boundary of $\Omega^{(r)}$ mapped forwards and backwards
by $\mathcal{J}_t$.

It is easy to see that the frontier stays on average near to the ``fluid limit'' in a different
coordinate system: $u=\ln(y)$ and $v=xy$. For this picture, see Figure \ref{fig:placeholder2}.

%

\begin{remark}[Connection to the Moran--Gamma subordinator]
The boundary sink point process $\Theta^{(r)}$ with scale-invariant intensity $r\,dy/y$ also provides a natural coupling across all balayage parameters $r \geq 0$. Specifically, let $\widetilde{\Theta}$ be a Poisson point process on $(0,\infty) \times (0,\infty)$ with intensity $(dy/y) \otimes dr$, and set $\Theta^{(r)} = \{y : (y, s) \in \widetilde{\Theta} \text{ for some } s \in [0, r]\}$. By Proposition~\ref{prop:sink_gamma}, the resulting horizontal shift $\mathsf{Z}^{(r)}$ carved out prior to the first surviving point in $\mathcal{R}$ evolves in $r$ as a standard Gamma subordinator, with 
$\mathsf{Z}^{(r)} \sim \mathrm{Gamma}(r,1)$ and independent increments. In this way, the scale-invariant boundary sinks realize the $\mathrm{Gamma}(r,1)$ shift via a Moran-type representation driven by a scale-invariant Poisson process (see Gnedin, Pitman, and Yor~\cite{GPY2006}, as well as Pitman~\cite{Pitman}), while the graphical flow $\mathcal{J}_t$ simultaneously preserves this structure in time $t$.
\end{remark}

\subsection{Pareto-record coupling and 1-balayage}
\label{sec:pareto-coupling}
Now we consider $r=R=1$.
Start with the unperturbed two-dimensional Poisson process $\Xi$ itself, without introducing independent randomness $\Theta^{(r)}$ or any sinks matched to points of $\Xi$, yet.
Instead, we consider a definition of the boundary sinks, based on $\Xi$ itself.
In other words, when conditioned on $\Xi$, our matching is conditionally deterministic.

Call a point $(x,y) \in \Xi$ a \emph{lower-left Pareto record} if
$$
 \Xi \cap \bigl((0,x] \times (0,y]\bigr) = \{(x,y)\}.
$$ 
Because $\Xi$ has unit intensity, these records form a bi-infinite chain $\{(x_k,y_k) : k \in \mathbb{Z}\}$ satisfying 
$$
\cdots < x_{-1} < x_0 < x_1 < \cdots\ \text{ and }\
\cdots > y_{-1} > y_0 > y_1 > \cdots\, .
$$
(We may choose $(x_0,y_0)$ so that $y_{-1}>1>y_0$.)
Conditioned on $(x_k,y_k)$, the next record $(x_{k+1},y_{k+1})$ is the leftmost point of $\Xi$ in the semi-infinite strip $(x_k,\infty) \times (0,y_k)$. 
So $y_{k+1} \sim \mathrm{Uniform}(0,y_k)$ and the log-gaps $(-\ln y_{k+1}) - (-\ln y_k)$ are i.i.d.\ $\mathrm{Exponential}(1)$ random variables. Since the hyperbolic flow $J_t(x,y) = (e^{-t}x, e^t y)$ preserves $\Xi$ and its Pareto records, projecting $(x_k,y_k) \mapsto (0,y_k)$ onto the vertical axis yields
\[
  \Theta^{(1)}_{\mathrm{Par}} := \{y_k : k \in \mathbb{Z}\},
\]
which is a stationary point process on $(0,\infty)$ with scale-invariant intensity $dy/y$---the exact marginal law of $\Theta^{(1)}$ as defined in Section~4.1 for the ``effective fractional $r$-balayage,'' because $r=1$ now.

Projecting onto the $y$-axis forgets the horizontal coordinates $(x_k)_{k \in \mathbb{Z}}$, whose conditional law given $\Theta^{(1)}_{\mathrm{Par}} = (y_k)_{k \in \mathbb{Z}}$ is governed by the same exponential shot-noise structure as Proposition~4.2: the increments $\Delta x_k := x_k - x_{k-1}$ are conditionally independent $\mathrm{Exponential}(y_{k-1})$ random variables, giving
\[
  x_k = \sum_{j=-\infty}^k \frac{\gamma_j}{y_{j-1}}, \qquad \gamma_j \stackrel{\text{i.i.d.}}{\sim} \mathrm{Exponential}(1).
\]
When $\Theta^{(1)}_{\mathrm{Par}}$ is paired with the full process $\Xi$, these dropped coordinates are recovered uniquely by the greedy matching algorithm of Section~4.1: because $(0,x_k] \times (0,y_k)$ is empty of points of $\Xi$, the leftmost point of $\Xi$ below height $\theta_k = y_k$ is precisely $(x_{k+1},y_{k+1})$. By induction, the matched set $\Omega_{\mathrm{det}}^{(1)}$ coincides identically with the set of all lower-left Pareto records of $\Xi$. 

Under $\mathcal{J}_t$, a point $(e^{-t}x_k, e^t y_k)$ lies in $\mathcal{R}$ and is matched to a sink above $y=1$ if and only if $y_k < e^{-t} < y_{k-1}$; at every $t \in \mathbb{R}$, this uniquely selects the instantaneous leader $\xi_1(t)$ of $\Pi(\mathcal{J}_t(\Xi) \cap \mathcal{R})$. Excavating $\Omega_{\mathrm{det}}^{(1)}$ therefore leaves the shifted sequence $(\xi_2(t), \xi_3(t), \dots)$ identically for all $t \in \mathbb{R}$.

There are two notable differences with the previous graphical representation.
First there is no auxiliary randomness. Instead $\Theta^{(1)}_{\mathrm{Par}}$ is determined
entirely from $\Xi$: they are coupled.
Secondly, the active count is $\mathfrak{M}_t \equiv 1$ rather than $\mathrm{Poisson}(1)$. 
One can think of the extra randomness in the previous picture, wherein $\mathfrak{M}_t \sim \mathrm{Poisson}(1)$, as compensating for the fact that $\Theta^{(r)}$ in the ``effective fractional balayage'' is not coupled to $\Xi$

One can also note that $\Omega^{(1)}_{\mathrm{Par}}$ is a decreasing sequence.
In the previous picture $\Omega^{(r)}$ is not strictly decreasing, not even for $r=1$. 
(But it is an $r$-decreasing sequence on average since $\mathbb{E}[\mathfrak{M}_t]=r$.)

Iterating this deterministic Pareto record extraction $R$ times on the remaining points $\Xi \setminus \bigcup_{j=1}^{R-1} \Omega_{\mathrm{Par}}^{(j)}$ yields $R$ coupled sink processes of intensity $dy/y$ whose greedy matchings excise $(\xi_1(t),\dots,\xi_R(t))$ at every time $t$, giving a deterministic graphical realization of integer $R$-point balayage.

\begin{figure}
    \centering
    \includegraphics[width=0.75\linewidth]{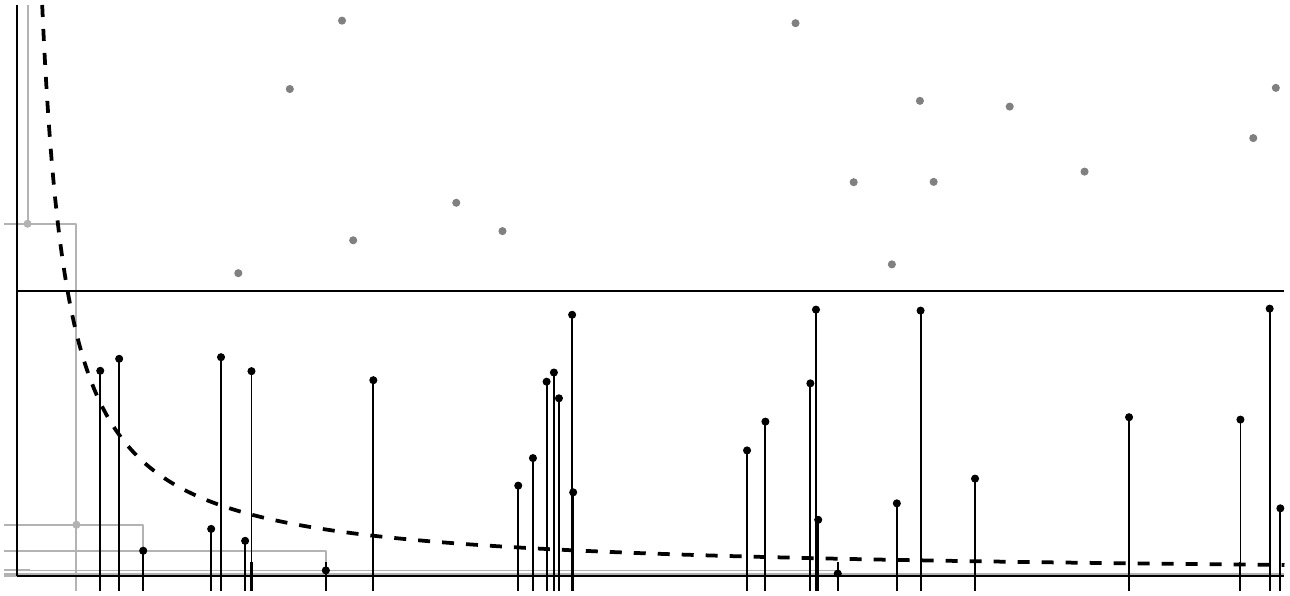}
    \caption{The graphical representation for the Pareto-coupling implementing 1-balayage:
    Now the ``boundary sinks'' are literally chosen from $\Xi$ itself using the Pareto records.
    This implements 1-balayage of the homogeneous Poisson point process of intensity 1. 
    Thus the resulting point process is the 2TTwALR$(2)$.}
    \label{fig:placeholder3}
\end{figure}

\subsection{Considerations for r in (-1,0)}
\label{sec:negative-r}
While the boundary-sink Poisson process $\Theta^{(r)}$ of Section~4.1 requires a non-negative intensity $r \geq 0$, the fractional balayage distribution $\mathrm{PD}^{(r)}(m,0)$ extends analytically to all $r > -1$ by replacing the shifted sequence $(Z^{(r)} + \xi_n)_{n \geq 1}$ in Definition~1.3 with $(\zeta_n^{(r)})_{n \geq 1}$, where $\zeta_1^{(r)} \sim \mathrm{Gamma}(1+r, 1)$ and the gaps $\zeta_{n+1}^{(r)} - \zeta_n^{(r)}$ ($n \geq 1$) are independent $\mathrm{Exponential}(1)$ random variables.

In the regime $-1 < r < 0$, classical Beta--Gamma algebra gives an exact static coupling to the unperturbed standard Poisson point process $(\xi_n)_{n \geq 1}$ on $(0,\infty)$: if $B \sim \mathrm{Beta}(1+r, -r)$ is independent of $(\xi_n)_{n \geq 1}$, then $\zeta_1^{(r)} := B\xi_1 \sim \mathrm{Gamma}(1+r, 1)$ and $W^{(r)} := (1-B)\xi_1 \sim \mathrm{Gamma}(-r, 1)$ are independent, so shifting every point leftward by the common offset $W^{(r)}$,
\[
  \zeta_n^{(r)} := \xi_n - (1-B)\xi_1 = B\xi_1 + (\xi_n - \xi_1), \qquad n \geq 1,
\]
realizes the $\mathrm{PD}^{(r)}(m,0)$ point process via a random contraction of the first gap rather than an additive Gamma excavation.

That all is at the level of a single-time marginal.
But one can also make a graphical representation for $-1<r<0$, albeit a different one
than for $r>0$.

At the process level under the hyperbolic flow $\mathcal{J}_t(x,y) = (e^{-t}x, e^t y)$ on $\mathcal{R} = (0,\infty) \times (0,1)$, the regime $-1 < r < 0$ is realized by \emph{probabilistic leader retention} at the exit boundary $y=1$ rather than boundary sinks on $x=0$. Specifically, evolve all points of $\Xi \cap \mathcal{R}$ under $\mathcal{J}_t$; whenever the instantaneous leftmost point $(x_1(t), y_1(t))$ hits the top boundary $y=1$ (which occurs at rate $1$), let it exit $\mathcal{R}$ with probability $1+r \in (0,1)$, and with complementary probability $-r \in (0,1)$ retain it at the same horizontal coordinate $x_1(t)$ while reassigning its vertical coordinate independently as $y_{\mathrm{new}} \sim \mathrm{Uniform}(0,1)$ (equivalently, restarting its \(\text{Exponential}(1)\) residual lifetime \(-\ln y_{\text{new}}\) at the same horizontal coordinate \(x_1(t)\)). Since retention preserves the horizontal position of the leader and leaves the unperturbed Poisson configuration on $(x_1(t), \infty) \times (0,1)$ untouched, the effective leader exit rate is reduced from $1$ to $1+r$, preserving the stationary $\mathrm{Gamma}(1+r, 1)$ law for $x_1(t)$ and hence the $\mathrm{PD}^{(r)}(m,0)$ partition structure for all $r \in (-1, 0)$.

To verify that all this is well-defined—without points accumulating near $x = 0$—and preserves both the $\mathrm{Gamma}(1+r,1)$ marginal of $x_1(t)$ and the full process-level $t$-stationarity, decorate the lead particle $(x_1(t), y_1(t)) \in \mathcal{R}$ with a non-negative integer counter $\alpha(t) \in \{0, 1, 2, \dots\}$ recording the number of times it has reached the exit boundary $y = 1$ and been re-randomized vertically into $(0,1)$. 

It is easy to see that, until it exits $\mathcal{R}$, the leader is not overtaken by any other
point of $\Xi \cap \mathcal{R}$ because $\mathcal{J}_t$ shrinks all $x$-coordinates by the factor
$e^{-t}$ uniformly, preserving the relative order of the $x$-coordinates.

Moreover, both a newly promoted leader (entering from the untouched Poisson bulk to the right of $x_1(t)$) and a retained leader (reassigned $y_{\mathrm{new}} \sim \mathrm{Uniform}(0,1)$) possess a vertical height distributed as $\mathrm{Uniform}(0,1)$, and hence an independent $\mathrm{Exponential}(1)$ residual lifetime $-\ln y$ until next hitting $y = 1$. Thus $\alpha_1(t)$ evolves as an autonomous continuous-time Markov chain on $\{0, 1, 2, \dots\}$ with transitions $k \to k+1$ at rate $-r$ (retention) and $k \to 0$ at rate $1+r$ (exit and replacement by a fresh leader with $\alpha = 0$). Writing $p_k(t) := \mathbb{P}(\alpha(t) = k)$ for $k \geq 0$, the Kolmogorov forward equations, starting from $\alpha_1(0) = 0$ ($p_0(0) = 1$), 
may be written in generating function format as
$$
G(z,t)\, \stackrel{\mathrm{def}}{:=}\, \sum_{k=0}^{\infty} p_k(t) z^k\, ,\ \text{ gives }\
\frac{\partial}{\partial t} G(z,t)\, 
=\, (1+r) - (1+rz) G(z,t)\, .
$$
This leads to 
$$
G(z,t)\, =\, \frac{1+r}{1+rz} - \frac{r(1-z)}{1+rz}\, e^{-(1+rz)t}\, .
$$
So the $(p_0(t),p_1(t),\dots)$ probability mass function
converges exponentially fast as $t \to \infty$ to the unique stationary $\mathrm{Geometric}(1+r)$ law: $\lim_{t \to \infty} p_k(t) = \pi_k = (1 + r)(-r)^k$ for $k=0,1,\dots$. 
Since this decays exponentially, we see that the negative moments $\mathbb{E}[x_1(t)^{-\beta}]$
are finite for some $\beta>0$ (See below: it follows as long as $\beta<1+r$.)

Finally, this geometric counter makes full $t$-stationarity transparent. Between exits, $x_1(t)$ drifts inward according to $dx_1/dt = -x_1$, while exits occur as a Poisson process of rate $1 + r > 0$ (the thinning of the rate-$1$ boundary-hitting clock by probability $1 + r$), at which $x_1(t)$ jumps rightward by $\Delta x \sim \mathrm{Exponential}(1)$. The infinitesimal generator of $x_1(t)$,
\[
\mathcal{A} f(x) = -x f'(x) + (1 + r)\int_0^{\infty} \bigl(f(x+u) - f(x)\bigr) e^{-u}\,du,
\]
has unique invariant density $x^r e^{-x}/\Gamma(1+r)$ on $(0,\infty)$, so $x_1(t) \sim \mathrm{Gamma}(1+r, 1)$ is stationary. Because $\mathcal{J}_t$ preserves $2\mathrm{d}$ Lebesgue measure on $(x_1(t), \infty) \times (0,1)$ and retention uses fresh independent vertical randomness, the conditional law of $\bigl(y_1(t), \Xi \cap ((x_1(t), \infty) \times (0,1))\bigr)$ given $x_1(t)$—namely $\mathrm{Uniform}(0,1)$ tensor a unit-intensity Poisson point process—is preserved for all $t \geq 0$.

\section{Summary and Outlook}
\label{sec:outlook}

In this article we have explored a special example of rank-dependent thinning called 2TTwALR$(1+r)$.
We are motivated by an economics application, where instead of companies interchanging their relative
positions, in our model some simply exit the market.
Then either extra regulatory attention or incumbency can cause the leader to exit the market
at a different rate than other companies.
One can imagine other applications, such as workers leaving a particular sector:
e.g., research mathematicians emigrating to work in industry in the AI sector.
Then the lead worker, as measured by research productivity, may receive extra encouragement to either leave or stay.

In these economic applications, we are obviously motivated by Fernholz’s “Atlas model” of interacting stock prices \cite{Fernholz} and the rank-dependent Brownian particle systems of Pal and Pitman \cite{PalPitman}. Unlike our discrete-event thinning dynamics—which permit a complete classification of time-stationary partition structures (Theorem 3.6) analogous to Aizenman and Ruzmaikina’s characterization of quasistationary measures for discrete-time competing particle systems \cite{AizenmanRuzmaikina}—characterizing the full family of invariant gap distributions and equilibrium fluctuations for the countably infinite continuous-time Atlas model is considerably more delicate due to local-time collisions, and was resolved only later by Dembo and Tsai \cite{DT17}, Sarantsev and Tsai \cite{ST17}, and Tsai \cite{Tsa18}.
Furthermore, while Dembo, Shkolnikov, Varadhan, and Zeitouni \cite{DSVZ16} established large deviation principles and hydrodynamic McKean–Vlasov limits for large systems of rank-interacting diffusions, Cabezas, Dembo, Sarantsev, and Sidoravicius \cite{CDSS19} analyzed the out-of-equilibrium convergence of the infinite Atlas model from general initial density profiles via one-sided Stefan free-boundary equations.

Our main result is that, for our own simple model of 2TTwALR$(1+r)$, if $R \in \{0,1,\dots\}$
then the solution is given by $R$-point balayage of the Poisson-Dirichlet $\mathrm{PD}(m,0)$ model,
by inspection. So it is natural to interpolate the formulas to obtain a solution
for general $r$, which we called $\mathrm{PD}^{(r)}(m,0)$.
An interesting feature is the existence of a Harris type graphical representation, which we have
explored.
Of course, because of the nature of thinning, there are ways to perturb one solution to obtain
another solution.
This was already noted in \cite{SVW}.
But this does not distract from the relative simplicity of the model.

We have considered dust-free random partition structures with infinitely many parts.
If there were finitely many parts, then eventually one would obtain an empty configuration
of points which is not normalizable.
Similarly, if there is any positive component of dust, then under thinning the dust comes to completely
dominate over time.

Finally, we note that several extensions and further developments point the way to open problems for further exploration.

\paragraph{A $\mathfrak{q}$-deformed thinning semigroup via $\mathfrak{q}$-Chu--Vandermonde.}
The proof of the semigroup property in Lemma~\ref{lem:injection} relies on the classical Chu--Vandermonde identity~\eqref{eq:ChuVandermonde}. We may consider replacing ordinary binomial coefficients and powers $(1 - e^{-t})^{k-1}$ by Gaussian $\mathfrak{q}$-binomial coefficients $\genfrac{[}{]}{0pt}{}{r+k-1}{k-1}_{\mathfrak{q}}$ and $\mathfrak{q}$-Pochhammer symbols $(e^{-t};\mathfrak{q})_{k-1} = \prod_{j=0}^{k-2}(1 - \mathfrak{q}^j e^{-t})$ for a deformation parameter $\mathfrak{q} \in (0,1)$. (We write 
this using Fraktur to avoid notational collision with the exponent in Section 1.)
This suggests defining a $\mathfrak{q}$-deformed negative binomial kernel
\[
\widetilde{Q}_{\mathfrak{q}}^{(r+1)}(k, t) := \genfrac{[}{]}{0pt}{}{r+k-1}{k-1}_{\mathfrak{q}} e^{-(r+1)t}\bigl(e^{-t};\mathfrak{q}\bigr)_{k-1}, \qquad k \in \mathbb{N}.
\]
By the $\mathfrak{q}$-Chu--Vandermonde identity (Corollary~10.2.2 of~\cite{AndrewsAskeyRoy}), these kernels satisfy a $\mathfrak{q}$-weighted convolution identity that defines a continuous convolution semigroup of random thinning injections $\pi_{\mathfrak{q}}^{(r)}(\cdot, t) \in \mathcal{F}^{\nearrow}(\mathbb{N}, \mathbb{N})$ in which deeper-ranked particles experience $\mathfrak{q}$-attenuated attrition rates rather than uniform unit-rate thinning. Characterizing the stationary random partition structures for this $\mathfrak{q}$-deformed thinning—and determining whether they relate to $\mathfrak{q}$-deformed Poisson--Dirichlet laws or $\mathfrak{q}$-regenerative compositions~\cite{GPY2006,GO}—is a natural direction for future work.

\paragraph{Multi-rank alternative thinning rates $(\lambda_1, \dots, \lambda_J)$.}
A  multi-parameter generalization of $\mathrm{2TTwALR}(1+r)$ assigns distinct thinning rates $\lambda_1, \dots, \lambda_J > 0$ to the top $J$ active ranks while keeping unit attrition rate $1$ for all subdominant ranks $n > J$. Under the rescaled thinning dynamics, the tail configuration $(\xi_{J+1}(t), \xi_{J+2}(t), \dots)$ conditioned on $\xi_J(t)$ remains an unperturbed unit-intensity Poisson point process on $(\xi_J(t), \infty)$. Consequently, the top $J$ positions $(x_1(t), \dots, x_J(t))$ evolve as an autonomous piecewise-deterministic Markov process on the chamber $0 < x_1 < \dots < x_J$ with inward drift $dx_i/dt = -x_i$ and jumps at rates $\lambda_\ell$ ($\ell = 1, \dots, J$) that delete $x_\ell$ and append $x_J + \Delta x$ ($\Delta x \sim \mathrm{Exponential}(1)$). In stationarity, the spacings $(x_1, x_2 - x_1, \dots, x_J - x_{J-1})$ admit an explicit rational Laplace transform—and when $\lambda_1 \geq \lambda_2 \geq \dots \geq \lambda_J \geq 1$, the process admits a tiered boundary-sink graphical representation with sink intensities $(\lambda_i - \lambda_{i+1})\,dy/y$ ($i = 1, \dots, J$, $\lambda_{J+1} := 1$)—yielding a $J$-parameter family of stationary laws $\mathrm{PD}^{(\lambda_1,\dots,\lambda_J)}(m,0)$ that extends $\mathrm{PD}^{(r)}(m,0)$ beyond rank-$1$ perturbations.

\paragraph{Econometrics.}
The most interesting future direction would be to try to link this model, or any
variant of it to the actual application we considered.
Namely,  to infer whether market share dynamics indicate that the leader faces additional regulatory scrutiny (\(r > 0\)) or enjoys network-effect incumbency (\(-1 < r < 0\)).
It is also possible that one and the same company experiences both at different times.
A primary goal is to explore mechanisms to connect the theoretical model to real-world
applications.

\section*{Acknowledgments}
We thank Google DeepMind's Gemini large language model for its utility as an exploratory tool during preliminary literature searches and initial algebraic calculations. In accordance with journal policy, all mathematical proofs, fluid limit derivations, and final prose were independently audited, verified, and digested by the human authors, who maintain full accountability for the integrity of this work.

We also thank Priya Mantraratnam for her work on an earlier version of this project.
\appendix



\baselineskip=12pt
\bibliographystyle{plain}

\begin{thebibliography}{10}





















\bibitem{Aldous2003}
D.~J.~Aldous.
\newblock A tractable complex network model based on the stochastic mean-field model of distance.
\newblock {\em Probability on Discrete Structures}, Encyclopaedia of Mathematical Sciences, 110:1--72, 2003.

\bibitem{AizenmanRuzmaikina}
M.~Aizenman and A.~Ruzmaikina.
\newblock Characterization of invariant measures at the leading edge for competing particle systems.
\newblock {\em Ann.~Probab.} {\bf 33}, no.~1, 82--113, 2005.

\bibitem{AndrewsAskeyRoy}
G.~E.~Andrews, R.~Askey, and R.~Roy.
\newblock {\em Special Functions: Encyclopecia of Mathematics and its Application, 71.}
\newblock Cambridge University Press, Cambridge, UK, 2000.

\bibitem{Billingsley}
P.~Billingsley.
\newblock {\em Probability and Measure: Third Edition.}
\newblock John Wiley \& Sons, New York, NY, 1995.

\bibitem{CDSS19}
M.~Cabezas, A.~Dembo, A.~Sarantsev, and V.~Sidoravicius. 
\newblock Brownian particles with rank-dependent drifts: out-of-equilibrium behavior. 
\newblock {\em Commun.~Pure Appl.~Math.} {\bf 72}, no.~7, 1424--1458, 2019.

\bibitem{DSVZ16}
A.~Dembo, M.~Shkolnikov, S.~R.~S.~Varadhan, and O.~Zeitouni. 
\newblock Large deviations for diffusions interacting through their ranks. 
\newblock {\em Commun.~Pure Appl.~Math.} {\bf 69}, no.~7, 1259--1313, 2016.

\bibitem{DT17}
A.~Dembo and L.~C.~Tsai.
\newblock Equilibrium fluctuation of the Atlas model.
\newblock {\em Ann.~Probab.} {\bf 45}, no.~6b, 4529--4560, 2017.

\bibitem{Dudley}
R.~M.~Dudley.
\newblock {\em Real Analysis and Probability.}
\newblock Cambridge University Press, Cambridge, UK, 2002.

\bibitem{Fernholz}
R.~Fernholz.
\newblock {\em Stochastic Portfolio Theory.}
\newblock Springer Verlag, New York, NY, 2002.

\bibitem{GO}
A.~Gnedin and G.~Olshanski.
\newblock $q$-exchangeability via quasi-invariance.
\newblock {\em Ann.~Probab.} {\bf 38}, no.~6, 2103--2135, 2010.

\bibitem{GPY2006}
A.~Gnedin, J.~Pitman, and M.~Yor.
\newblock Asymptotic laws for regenerative compositions: gamma subordinators and the like.
\newblock {\em Probability Theory and Related Fields}, 135(4):576--602, 2006.


\bibitem{Harris}
T.~E.~Harris
\newblock Additive set-valued Markov processes and graphical methods.
\newblock {\em Ann.~Probab.} {\bf 6}, no.~3, 355--378 (1978)

 
\bibitem{Kallenberg}
O.~Kallenberg.
\newblock {\em Random measures.}
\newblock Academic Press, Inc., London, 1983.

\bibitem{Liggett}
T.~Liggett.
\newblock {\em Interacting Particle Systems. Grundlehren der mathematischen Wissenschaften 276.}
\newblock Springer-Verlag, New York, 1985. 

\bibitem{LOS}
J.~Lindenstrauss, G.~Olsen and Y.~Sternfeld.
\newblock The Poulsen simplex.
\newblock {\em Ann.~Inst.~Fourier (Grenoble)} {\bf 28}, no.~1,
91--114 (1978).


\bibitem{PalPitman}
S.~Pal and J.~Pitman.
\newblock One-dimensional Brownian particle systems with rank-dependent drifts.
\newblock {\em Ann.~Appl.~Probab.} {\bf 18}, no.~6, 2179--2207, 2008.

\bibitem{Pitman}
J.~Pitman.
\newblock {\em Combinatorial Stochastic Processes, \'{E}cole d'Et\'{e} de Probabilit\'{es} de Saint-Flour
XXXII-2002.}
\newblock Springer-Verlag, Berlin 2006.

\bibitem{PitmanYor}
J.~Pitman and M.~Yor.
\newblock The two-parameter Poisson-Dirichlet distribution
derived from a stable subordinator.
\newblock {\em Ann.~Probab.} {\bf 25}, no.~2, 855--900 (1997).

\bibitem{Ruelle}
D.~Ruelle.
A mathematical reformulation of Derrida's REM and GREM.
\newblock {\em Commun.~Math.~Phys.} {\bf 108}, 225--239 (1987).

\bibitem{ST17}
A.~Sarantzev and L.-C.~Tsai.
\newblock  Stationary gap distributions for infinite systems of competing Brownian particles.
\newblock {\em Electron.~J.~Probab.} {\bf 22}, paper no.~56, 1--20, 2017.

\bibitem{SVW}
S.~Starr, B.~Vermesi and A.~Wei.
\newblock About thinning invariant partition structures.
\newblock {\em J.~Statist.~Phys.} {\bf 148}, no.~2, 325--344 (2012).

\bibitem{Tsa18}
L.-C.~Tsai.
\newblock Stationary distributions of the Atlas model. 
\newblock {\em Electron.~Commun.~Probab.} {\bf 23}, paper no.~10, 1–10, 2018.





\end{thebibliography}

\appendix

\section{Matlab m-file}

Regular graphical representation in $(0,\infty) \times (0,\infty)$.
{
\tiny

\begin{verbatim}
n=1000;
N=n*n;
xwin=25;
ywin=2;
xblk=[];
yblk=[];
K = poissrnd(N);
x_lst = n*rand(K,1);
y_lst = n*rand(K,1);
r = 3;          % Set fractional balayage parameter r directly
L = 100;            % Cutoff in log-height: y in (exp(-L), 1)
bdry_num = poissrnd(r * L);
bdry_logs = -L * rand(bdry_num, 1);
bdry_pts=n*exp(bdry_logs);
figure
hold on
axis([0,1,exp(-5),1])
sort_pts = sort(bdry_pts,'descend');
% RGB triplet for a medium-light gray (e.g., [0.5 0.5 0.5])
grayColor = [0.5 0.5 0.5];
for k=1:length(sort_pts)
    t = sort_pts(k);
    Idx = find(y_lst<t);
    if length(Idx) == 0
        plot([0,0],[t,t],'Color',grayColor)
    else
        x_lst_select = x_lst(Idx);
        [x,idx] = min(x_lst_select);
        K = Idx(idx);
        y = y_lst(K);
        plot([0,x,x],[t,t,y],'Color',grayColor)
        if t<1
            xblk = [xblk,x];
            yblk = [yblk,y];
        else
            plot(x,y,'.','Color',grayColor,MarkerSize=10)
            if y<1
                plot([x,x],[y,0],'Color',grayColor)
                plot([x,x],[-0.05,0.05],'Color',grayColor,LineWidth=1)
            end
        end
        x_lst(K)=[];
        y_lst(K)=[];
    end
end
for j=1:length(bdry_pts)
    t = bdry_pts(j);
    plot([-0.25,0.25],[t,t],'Color',grayColor,LineWidth=1)
end
for j=1:length(xblk)
    x=xblk(j);
    y=yblk(j);
    if x<xwin && y<ywin
        if y<1
            plot(x,y,'k.',MarkerSize=10)
            plot([x,x],[y,0],'k')
            plot([x,x],[-0.05,0.05],'k',LineWidth=1)
        else
            plot(x,y,'.','Color',grayColor,MarkerSize=10)
        end
    end
end
for j=1:length(x_lst)
    x=x_lst(j);
    y=y_lst(j);
    if x<xwin && y<ywin
        if y<1
            plot(x,y,'k.',MarkerSize=10)
            plot([x,x],[y,0],'k')
            plot([x,x],[-0.05,0.05],'k',LineWidth=1)

        else
            plot(x,y,'.','Color',grayColor,MarkerSize=10)
        end
    end
end
xstep=0:0.001:100;
plot(xstep,r./xstep,'k--',LineWidth=2)
plot([0,n],[1,1],'k')
axis([-0.5,25,-0.2,2])
plot([0,xwin],[0,0],'k')
plot([0,0],[0,ywin],'k')
axis off

\end{verbatim}
}

Hyperbolic coordinate graphical representation (by Google DeepMind Gemini)

{\tiny

\begin{verbatim}
% Stationary matching and transformed Gamma envelope in (u, v) coordinates
% u = ln(y) (horizontal axis), v = x*y (vertical axis)

r = 3;                  % Fractional balayage parameter r > 0
u_win = [-5, 5];        % Visible horizontal window for u = ln(y)
v_win = [0, r + 4*sqrt(r) + 5]; % Visible vertical window for v = x*y

% Buffers to the left and right of u_win so edge effects vanish inside u_win
u_min = u_win(1) - 8;
u_max = u_win(2) + 8;
v_max = v_win(2) + 4*sqrt(r) + 8;

% 1. Generate homogeneous Poisson processes directly in (u, v)
N_bulk = poissrnd(v_max * (u_max - u_min));
v_bulk = v_max * rand(N_bulk, 1);
u_bulk = u_min + (u_max - u_min) * rand(N_bulk, 1);

[u_bulk, sort_idx] = sort(u_bulk, 'descend');
v_bulk = v_bulk(sort_idx);
logv_bulk = log(v_bulk);
matched = false(N_bulk, 1);

N_sink = poissrnd(r * (u_max - u_min));
u_sink = sort(u_min + (u_max - u_min) * rand(N_sink, 1), 'descend');

sink_u = nan(N_sink, 1);
bulk_v = nan(N_sink, 1);
bulk_u = nan(N_sink, 1);
corner_x = nan(N_sink, 1); % Original x-coordinate of each elbow corner

% 2. Greedy matching minimizing x = v * exp(-u) via log(v) - u
search_depth = 10;
ptr_top = 1;

for k = 1:N_sink
    u0 = u_sink(k);
    while ptr_top <= N_bulk && u_bulk(ptr_top) >= u0
        ptr_top = ptr_top + 1;
    end
    if ptr_top > N_bulk, break; end
    
    ptr_bot = ptr_top;
    while ptr_bot <= N_bulk && u_bulk(ptr_bot) > u0 - search_depth
        ptr_bot = ptr_bot + 1;
    end
    
    cand = ptr_top:(ptr_bot - 1);
    cand = cand(~matched(cand));
    
    if ~isempty(cand)
        [~, best_local] = min(logv_bulk(cand) - u_bulk(cand));
        best_idx = cand(best_local);
        
        matched(best_idx) = true;
        sink_u(k) = u0;
        bulk_v(k) = v_bulk(best_idx);
        bulk_u(k) = u_bulk(best_idx);
        corner_x(k) = v_bulk(best_idx) * exp(-u_bulk(best_idx));
    end
end

% 3. Plot sinks, matched bulk points, and elbow connectors in (u, v)
figure('Color', 'w');
hold on;

grayColor = [0.65, 0.65, 0.65];
u_arc = linspace(0, 1, 30);

for k = 1:N_sink
    u0 = sink_u(k);
    if isnan(u0), continue; end
    v1 = bulk_v(k);
    u1 = bulk_u(k);
    
    if u1 > u_win(2) || u0 < u_win(1)
        continue;
    end
    
    v_corner = corner_x(k) * exp(u0);
    
    % Vertical segment from boundary v=0 up to elbow corner (u0, v_corner)
    plot([u0, u0], [0, v_corner], 'Color', grayColor, 'LineWidth', 0.9);
    
    % Exponential arc from (u0, v_corner) to matched bulk point (u1, v1)
    u_seg = u1 + (u0 - u1) * u_arc;
    v_seg = v1 * exp(u_seg - u1);
    plot(u_seg, v_seg, 'Color', grayColor, 'LineWidth', 0.9);
    
    % Sink tick on v = 0 and matched bulk point
    plot([u0, u0], [-0.15, 0.15], 'k-', 'LineWidth', 1);
    plot(u1, v1, 'k.', 'MarkerSize', 11);
end

% 4. Construct and plot the transformed Gamma envelope
valid = ~isnan(sink_u);
u_s = sink_u(valid);
x_c = corner_x(valid);
x_env = cummax(x_c);

for k = 1:(length(u_s) - 1)
    u_top = u_s(k);
    u_bot = u_s(k+1);
    
    if u_bot > u_win(2) || u_top < u_win(1)
        continue;
    end
    
    % Exponential curve v = x_env(k) * exp(u) between u_bot and u_top
    u_seg = linspace(u_bot, u_top, 25);
    v_seg = x_env(k) * exp(u_seg);
    plot(u_seg, v_seg, 'k-', 'LineWidth', 1.8);
    
    % Vertical jump at u_bot if the next corner extends higher in v
    if x_env(k+1) > x_env(k)
        v_low = x_env(k) * exp(u_bot);
        v_high = x_env(k+1) * exp(u_bot);
        plot([u_bot, u_bot], [v_low, v_high], 'k-', 'LineWidth', 1.8);
    end
end

% Horizontal boundary v = 0 and hydrodynamic line v = r
plot(u_win, [0 0], 'k-', 'LineWidth', 1.2);
plot(u_win, [r, r], 'k--', 'LineWidth', 1.8);

xlim(u_win);
ylim([-0.5, v_win(2)]);
xlabel('u = ln(y)');
ylabel('v = x y');
box off;
\end{verbatim}

}

\end{document}